\documentclass[12pt,letterpaper,reqno]{amsart}

\usepackage{times}
\usepackage[T1]{fontenc}
\usepackage{mathrsfs}
\usepackage{latexsym}
\usepackage[dvips]{graphics}
\usepackage{epsfig}
\usepackage{amsmath,amsfonts,amsthm,amssymb,amscd}
\input amssym.def
\input amssym.tex

\newcommand{\hphipr}{\hphi\left(\frac{\log p}{\log R}\right)}

\newcommand{\hkn}{H_k^\ast(N)}

\newcommand{\hphi}{\widehat{\phi}}  

\newcommand\be{\begin{equation}}
\newcommand\ee{\end{equation}}
\newcommand\bea{\begin{eqnarray}}
\newcommand\eea{\end{eqnarray}}
\newcommand\bi{\begin{itemize}}
\newcommand\ei{\end{itemize}}
\newcommand\ben{\begin{enumerate}}
\newcommand\een{\end{enumerate}}
\newcommand\bc{\begin{center}}
\newcommand\ec{\end{center}}
\newcommand\ba{\begin{array}}
\newcommand\ea{\end{array}}

\def\notdiv{\ \mathbin{\mkern-8mu|\!\!\!\smallsetminus}}

\newcommand{\R}{\ensuremath{\mathbb{R}}}

\newcommand{\Q}{\mathbb{Q}}

\newcommand{\F}{\mathcal{F}}

\newcommand{\foh}{\frac{1}{2}}  

\newcommand{\es}[2]{\left\langle {#1 \atop #2} \right\rangle}

\newcommand{\js}[1]{{\underline{#1}\choose p}}

\newtheorem{thm}{Theorem}[section]

\newtheorem{lem}[thm]{Lemma}

\newtheorem{defi}[thm]{Definition}

\theoremstyle{definition}

\newcommand{\ncr}[2]{{#1 \choose #2}}
\newcommand{\twocase}[5]{#1 \begin{cases} #2 & \text{{\rm #3}}\\ #4
&\text{{\rm #5}} \end{cases}   }

\newcommand{\gep}{\epsilon}

\newcommand{\fwf}{\frac1{W_R(\F)}}

\newcommand{\glfp}{\lambda_f(p)}

\newcommand{\gafp}{\alpha_f(p)}

\newcommand{\gbfp}{\beta_f(p)}
\newcommand{\glf}{\lambda_f}
\newcommand{\lp}{\lambda_p}

\newcommand{\lix}[1]{{\rm Li}_{#1}(x)}

\begin{document}

\title{An identity for sums of polylogarithm functions}

\author{Steven J. Miller}
\email{sjmiller@math.brown.edu} \address{Department of Mathematics,
Brown University, Providence, RI 02912} 

\date{\today}

\thanks{The author would like to thank Walter Becker and Eduardo
Due\~nez for useful discussions, Toufik Mansour for catching a typo in an earlier draft, and his son Cam and nephew Eli Krantz for sleeping quietly on his arm while some of the calculations were performed. Many of the formulas for expressions in this paper were first guessed by using Sloane's On-Line Encyclopedia of Integer
Sequences \cite{Sl}. The author was partly supported by NSF grant
DMS0600848.}

\begin{abstract} We derive an identity for certain linear
combinations of polylogarithm functions with negative exponents,
which implies relations for linear combinations of Eulerian numbers.
The coefficients of our linear combinations are related to expanding
moments of Satake parameters of holomorphic cuspidal newforms in
terms of the moments of the corresponding Fourier coefficients, which has applications in analyzing lower order terms in the behavior of zeros of $L$-functions near the central point.
\end{abstract}

\maketitle



\section{Introduction}
\setcounter{equation}{0}

The polylogarithm function $\lix{s}$ is \be \lix{s} \ = \
\sum_{k=1}^\infty k^{-s} x^k. \ee If $s$ is a negative integer, say
$s=-r$, then the polylogarithm function converges for $|x|<1$ and
equals \be\label{eq:eulerinpolylog} \lix{-r} \ = \
\frac{\sum_{j=0}^r \es{r}{j} x^{r-j}}{(1-x)^{r+1}},\ee where the
$\es{r}{j}$ are the Eulerian numbers. The Eulerian number
$\es{r}{j}$ is the number of permutations of $\{1,\dots,r\}$ with
$j$ permutation ascents. One has \be \es{r}{j} \ = \
\sum_{\ell=0}^{j+1} (-1)^\ell \ncr{r+1}{\ell}(j-\ell+1)^r. \ee We
record $\lix{-r}$ for some $r$:\bea \lix{0} & \ = \ & \frac{x}{1-x}
\nonumber\\ \lix{-1} & = & \frac{x}{(1-x)^2} \nonumber\\
\lix{-2} & = & \frac{x^2+x}{(1-x)^3} \nonumber\\
\lix{-3} & = & \frac{x^3+4x^2+x}{(1-x)^4}\nonumber\\
\lix{-4} & = & \frac{x^4+11x^3+11x^2+x}{(1-x)^5} \nonumber\\
\lix{-5} & = & \frac{x^5+26x^4+66x^3+26x^2+x}{(1-x)^6}
.\eea

From \eqref{eq:eulerinpolylog} we immediately deduce that, when $s$
is a negative integer, $\lix{s}$ is a rational function whose
denominator is $(1-x)^{|s|}$. Thus an appropriate integer linear
combination of $\lix{0}$ through $\lix{-n}$ should be a simple
rational function. In particular, we prove

\begin{thm}\label{thm:lincombpolylog} Let $a_{\ell,i}$ be the coefficient of
$k^{i}$ in $\prod_{j=0}^{\ell-1} (k^2-j^2)$, and let $b_{\ell,i}$ be
the coefficient of $k^{i}$ in $(2k+1)\prod_{j=0}^{\ell-1}
(k-j)(k+1+j)$. Then for $|x| < 1$ and $\ell \ge 1$ we have \bea
a_{\ell,2\ell} \lix{-2\ell} + \cdots + a_{\ell,0}\lix{0} &\ =
\ &  \frac{(2\ell)!}{2}\ \frac{x^\ell(1+x)}{(1-x)^{2\ell+1}} \nonumber\\
b_{\ell,2\ell+1} \lix{-2\ell-1} + \cdots + b_{\ell,0}\lix{0}& \ = \
& (2\ell+1)!\ \frac{x^\ell(1+x)}{(1-x)^{2\ell+2}}. \eea
\end{thm}

We prove Theorem \ref{thm:lincombpolylog} in
\S\ref{sec:valuesBrpupto6}. While Theorem \ref{thm:lincombpolylog}
only applies to linear combinations of polylogarithm functions with
$s$ a negative integer, it is interesting to see how certain special
combinations equal a very simple rational function. One application
is to use this result to deduce relations among the Eulerian numbers
(possibly by replacing $x$ with $1-x$ when expanding); another is of
course to write $\lix{-n}$ in terms of $\lix{-n+1}$ through
$\lix{0}$. The coefficients $a_{\ell,i}$ and $b_{\ell,i}$ which
occur in our linear combinations also arise in expressions involving
the Fourier coefficients of cuspidal newforms. We describe
this connection in greater detail in \S\ref{sec:connwithNT}; these expansions are related to understanding the lower order terms in the behavior of zeros of $L$-functions of cuspidal newforms near the central point. (see \cite{Mil3} for a complete analysis).


\section{Proof of Theorem \ref{thm:lincombpolylog}}\label{sec:valuesBrpupto6} \setcounter{equation}{0}

Before proving Theorem \ref{thm:lincombpolylog} we introduce some
useful expressions.

\begin{defi}\label{def:formscmr} Let \bea c_{2\ell}
 \ = \  \prod_{j=0}^{\ell-1} (\ell^2 - j^2)\ = \ (2\ell)! / 2, \ \
\ c_{2\ell+1} \ = \ (2\ell+1)\prod_{j=0}^{\ell-1} (\ell-j)(\ell+1+j)
\ = \ (2\ell+1)!.\eea Define constants $c_{m,r}$ as follows:
$c_{m,r} = 0$ if $m \not\equiv r \bmod 2$, and \ben
\item for $r$ even, $c_{0,0} = 0$, $c_{2k,0} = (-1)^k 2$ for
$k\ge 1$, and for $1 \le \ell \le k$ set \be c_{2k,2\ell} \ = \
\frac{(-1)^{k+\ell}}{c_{2\ell}}\prod_{j=0}^{\ell-1} (k^2 - j^2) \ =
\ \frac{(-1)^{k+\ell}}{c_{2\ell}}\frac{k \cdot
(k+\ell-1)!}{(k-\ell)!}; \ee
\item for $r$ odd and $0 \le \ell \le k$ set \bea c_{2k+1,2\ell+1} \ = \
\frac{(-1)^{k+\ell}}{c_{2\ell+1}}(2k+1) \prod_{j=0}^{\ell-1}
(k-j)(k+1+j) \ = \ \frac{(-1)^{k+\ell}(2k+1)}{c_{2\ell+1}}
\frac{(k+\ell)!}{(k-\ell)!}. \eea \een Note $c_{m,r}=0$ if $m<r$.
Finally, set $B_r(x) = \sum_{m=0}^\infty c_{m,r} (-x)^{m/2}$ for
$|x|<1$. Thus for $r = 2\ell\ge 2$ we have \be B_{2\ell}(x) \ = \
\sum_{m = 0}^\infty c_{m,2\ell} (-x)^{m/2} \ = \ \sum_{k=1}^\infty
\left( \frac{(-1)^{k+\ell}}{c_{2\ell}} \prod_{j=0}^{\ell-1}
(k^2-j^2)\right) (-x)^k.\ee
 \end{defi}

Immediately from the definition of $c_{r}$ we have
\be\label{eq:c2lmo2l} c_{2\ell-1} \ = \ \frac{c_{2\ell}}{\ell}\ = \
\frac{c_{2\ell+1}}{2\ell(2\ell+1)}, \ee as well as
\be\label{eq:crelsconjpr} c_{2\ell+2} \ = \
(2\ell+2)(2\ell+1)c_{2\ell}, \ \ \ c_{2\ell+3} \ = \
(2\ell+3)(2\ell+2)c_{2\ell+1}. \ee

While the definition of the $c_{m,r}$'s above may seem arbitrary,
these expressions arise in a very natural manner in number theory.
See \cite{Mil3} for applications of these coefficients in
understanding the behavior of zeros of ${\rm GL}(2)$ $L$-functions;
we briefly discuss some of these relations in
\S\ref{sec:connwithNT}.

\begin{proof}[Proof of Theorem \ref{thm:lincombpolylog}]
We first consider the case of $r=2\ell$ even. We proceed by
induction. We claim that \bea\label{eq:claimtoprove} B_{2\ell}(x) &
\ = \ & \sum_{k=1}^\infty \left( \frac{(-1)^{k+\ell}}{c_{2\ell}}
\prod_{j=0}^{\ell-1} (k^2-j^2)\right) (-x)^k \ = \ (-1)^\ell
\frac{x^\ell(1+x)}{(1-x)^{2\ell+1}} \eea for all $\ell$.

We consider the basis case, when $\ell = 1$. Thus we must show for
$|x| < 1$ that $B_2(x) = -x(1+x)/(1-x)^3$. As $r=2$, the only
non-zero terms are when $m=2k > 0$ is even. As $c_2 = 2$ and
$c_{2k,2} = (-1)^{k+1} k^2$ for $k\ge 1$, we find that \bea B_2(x) &
\ = \ & \sum_{k=1}^\infty (-1)^{k+1} k^2 (-x)^{2k/2} \ = \ -
\sum_{k=1}^\infty k^2 x^k \ =\ -\lix{-2} \ = \
-\frac{x(1+x)}{(1-x)^3},\eea which completes the proof of the basis
step. For the inductive step, we assume \be\label{eq:evenfirststep}
\sum_{k=1}^\infty \left( \frac{(-1)^{k+\ell}}{c_{2\ell}}
\prod_{j=0}^{\ell-1} (k^2-j^2)\right) (-x)^k \ = \ (-1)^\ell
\frac{x^\ell(1+x)}{(1-x)^{2\ell+1}}, \ee and we must show the above
holds with $\ell$ replaced by $\ell+1$. We apply the differential
operator \be \left(x \frac{d}{dx}\right)^2 - \ell^2 \ee to both
sides of \eqref{eq:evenfirststep}. After canceling the minus signs
we obtain \bea \sum_{k=1}^\infty \left( c_{2\ell}^{-1}
\prod_{j=0}^{\ell-1} (k^2-j^2)\right) (k^2-\ell^2) x^k & \ = \ &
\left(\left(x \frac{d}{dx}\right)^2-\ell^2\right)\left(
\frac{x^\ell(1+x)}{(1-x)^{2\ell+1}}\right)\nonumber\\
\sum_{k=1}^\infty c_{2\ell}^{-1}\left(\prod_{j=0}^{\ell}
(k^2-j^2)\right) x^k & \ = \ & (2\ell+2)(2\ell+1)
\frac{x^{\ell+1}(1+x)}{(1-x)^{2(\ell+1)+1}} \nonumber\\
\sum_{k=1}^\infty c_{2(\ell+1)}^{-1}\left(\prod_{j=0}^{\ell+1-1}
(k^2-j^2)\right) x^k & \ = \ &
\frac{x^{\ell+1}(1+x)}{(1-x)^{2(\ell+1)+1}}, \eea where the last
line follows from \eqref{eq:crelsconjpr}, which says $c_{2\ell+2}
 = (2\ell+2)(2\ell+1)c_{2\ell}$. Thus \eqref{eq:claimtoprove} is true for
all $\ell$.

As we have defined $a_{\ell,i}$ to be the coefficient of $k^{i}$ in
$\prod_{j=0}^{\ell-1} (k^2-j^2)$, \eqref{eq:claimtoprove} becomes
\bea \sum_{k=1}^\infty \sum_{i = 0}^{2\ell} a_{\ell,i}\ k^{i}\ x^k &
\ = \ & c_{2\ell} \frac{x^{\ell}(1+x)}{(1-x)^{2\ell+1}}.
\eea The proof of Theorem \ref{thm:lincombpolylog} for $r$ even is
completed by noting that the left hand side above is just \be
a_{\ell,2\ell} \lix{-2\ell} + \cdots + a_{\ell,0}\lix{0}. \ee

The proof for $r=2\ell+1$ odd proceeds similarly, the only
significant difference is that now we apply the operator
\be\label{eq:diffop} \left(x\frac{d}{dx}\right)^2 \ + \
\left(x\frac{d}{dx}\right) \ - \ \ell(\ell+1), \ee which will bring
down a factor of $(k-\ell)(k+1-\ell)$.
\end{proof}

\section{Connections with number theory}\label{sec:connwithNT}
\setcounter{equation}{0}

We now describe how our polylogarithm identity can be used to analyze zeros of $L$-functions near the central point. Katz and Sarnak \cite{KaSa} conjecture that, in the limit as
the conductors tend to infinity, the behavior of the normalized
zeros near the central point agree with the $N\to\infty$ scaling
limit of the normalized eigenvalues near $1$ of a subgroup of
$U(N)$ ($N\times N$ unitary matrices); see \cite{DM,FI,Gu,HR,HM,ILS,KaSa,Mil1,Ro,Rub,Yo} for many examples. While the main terms for many families are the same as the
conductors tend to infinity, a more careful analysis of the explicit
formula allows us to isolate family dependent lower order terms.

Our coefficients $c_{m,r}$ are related to writing the moments of
Satake parameters of certain ${\rm GL}(2)$ $L$-functions in terms of
the moments of their Fourier coefficients, which we briefly review.
Let $H^\star_k(N)$ be the set of all holomorphic cuspidal newforms
of weight $k$ and level $N$; see \cite{Iw2} for more details. Each
$f\in H^\star_k(N)$ has a Fourier expansion
\begin{equation}\label{eq:defLsf1}
f(z)\ =\ \sum_{n=1}^\infty a_f(n) e(nz).
\end{equation}
Let $\lambda_f(n) =  a_f(n) n^{-(k-1)/2}$. These coefficients
satisfy multiplicative relations, and $|\lambda_f(p)| \le 2$. The
$L$-function associated to $f$ is
\begin{equation}
L(s,f)\ =\ \sum_{n=1}^\infty \frac{\lambda_f(n)}{n^{s}} \ = \
\prod_p \left(1 - \frac{\lambda_f(p)}{p^s} +
\frac{\chi_0(p)}{p^{2s}}\right)^{-1},
\end{equation} where $\chi_0$ is the principal character with
modulus $N$. We write \be \lambda_f(p) \ = \ \alpha_f(p) +
\beta_f(p). \ee For $p\ \notdiv N$, $\alpha_f(p)\beta_f(p) = 1$ and
$|\alpha_f(p)| = 1$. If $p|N$ we take $\alpha_f(p) = \lambda_f(p)$
and $\beta_f(p) = 0$. Letting \bea L_\infty(s,f) & \ = \ &
\left(\frac{2^k}{8\pi}\right)^{1/2}\
\left(\frac{\sqrt{N}}{\pi}\right)^s\
\Gamma\left(\frac{s}2+\frac{k-1}4\right)\
\Gamma\left(\frac{s}2+\frac{k+1}4\right) \eea denote the local
factor at infinity, the completed $L$-function is
\begin{equation}\label{eq:defLsf2}
\Lambda(s,f) \ =\ L_\infty(s) L(s,f) \ =\ \epsilon_f \Lambda(1-s,f),
\ \ \ \epsilon_f = \pm 1.\end{equation}

The zeros of $L$-functions often encode arithmetic information, and
their behavior is well-modeled by random matrix theory
\cite{CFKRS,KaSa,KeSn3}. The main tool in analyzing the behavior of these
zeros is through an explicit formula, which relates sums of a test
function at these zeros to sums of the Fourier transform of the test
function at the primes, weighted by factors such as $\alpha_f(p)^m +
\beta_f(p)^m$. For example, if $\phi$ is an even Schwartz function,
$\hphi$ its Fourier transform, and $\foh + i\gamma_f$ denotes a
typical zero of $\Lambda(s,f)$ for $f\in H^\star_k(N)$ (the
Generalized Riemann Hypothesis asserts each $\gamma_f \in \R$), then
the explicit formula is \bea\label{eq:explicitformula} &  & \frac1{|H^\ast_k(N)|} \sum_{f\in
H^\ast_k(N)} \sum_{\gamma_f} \phi\left(\gamma_f \frac{\log
N}{2\pi}\right)  \nonumber\\ & & \ \ \ \ \ = \ \frac{A(\phi)}{\log
N} \ +\ \frac1{|H^\ast_k(N)|} \sum_{f\in H^\ast_k(N)}
\sum_{m=1}^\infty \sum_p \frac{\alpha_f(p)^m +
\beta_f(p)^m}{p^{m/2}} \frac{\log p}{\log N}\ \hphi\left(m\frac{\log
p}{\log N}\right); \eea see \cite{ILS,Mil3} for details and a definition
of $A(\phi)$. Similar expansions hold for other families of
$L$-functions. Information about the distribution of zeros in a
family of $L$-functions (the left hand side above) is obtained by
analyzing the prime sums weighted by the moments of the Satake
parameters (on the right hand side). Thus it is important to be able
to evaluate quantities such as \be \frac1{|\mathcal{F}|} \sum_{f\in
\mathcal{F}} \left(\alpha_f(p)^m + \beta_f(p)^m\right) \ee for
various families of $L$-functions.

For some problems it is convenient to rewrite
$\alpha_f(p)^m + \beta_f(p)^m$ in terms of a polynomial in $\lambda_f(p)$. This
replaces moments of the Satake parameters $\alpha_f(p)$ and
$\beta_f(p)$ with moments of the Fourier coefficients
$\lambda_f(p)$, and for many problems the Fourier coefficients are
more tractable; we give two examples.

First, the $p$\textsuperscript{th} coefficient of the $L$-function of the elliptic curve $y^2 = x^3 + Ax + B$ is $p^{-1/2}$ $\sum_{x \bmod p}$ $\js{x^3+Ax+B}$; here $\js{x}$ is the Legendre symbol, which is 1 if $x$ is a non-zero square modulo $p$, $0$ if $x \equiv 0 \bmod p$, and $-1$ otherwise. Our sum equals the number of solutions to $y^2 \equiv x^3 + Ax + B \bmod p$, and thus
these sums can be analyzed by using results on sums of Legendre symbols (see for example \cite{ALM,Mil2}).

Second, the Petersson formula
(see Corollary 2.10, Equation (2.58) of \cite{ILS}) yields, for $m,
n > 1$ relatively prime to the level $N$,
\be\label{eq:PeterssonFormula} \fwf\sum_{f \in \hkn} w_R(f)
\glf(m)\glf(n) \ = \ \delta_{mn} \ + \ O\left((mn)^{1/4}\frac{\log
2mnN}{k^{5/6} N}\right),\ee where $\delta_{mn} = 1$ if $m=n$ and $0$
otherwise. Here the $w_R(f)$ are the harmonic weights \be w_R(f) \ = \ \zeta_N(2) /
Z(1,f) \ = \ \zeta(2) / L(1,{\rm sym}^2 f). \ee They are mildly varying, with (see \cite{Iw1,HL})  \be N^{-1-\gep}\ \ll_k \ \omega_R(f) \ \ll_k \
N^{-1+\gep}; \ee if we allow ineffective constants we can replace
$N^\gep$ with $\log N$ for $N$ large.

We can now see why our polylogarithm identity is useful.
Using $\gafp + \gbfp = \lp$, $\gafp \gbfp = 1$ and $|\gafp| =
|\gbfp| = 1$, we find that \bea
                    \begin{array}{rrlrlrlrlrr}
\gafp^{\  } + \gbfp^{\  }    & = & \glfp &  &  &  &  &  &  &  &  \\ \ \\
\gafp^2+\gbfp^2 & = & \glfp^2 & - & 2 &  &  &  &  &  &  \\ \ \\
\gafp^3+\gbfp^3 & = & \glfp^3 & - & 3\glfp  &  &  &  &  &  &  \\ \ \\
\gafp^4+\gbfp^4 & = & \glfp^4 & - & 4\glfp^2 & + & \ \  2 &  &  &  &  \\
\ \\
\gafp^5+\gbfp^5 & = & \glfp^5 & - & 5\glfp^3 & + & \ \ 5\glfp &  & &
&
\\ \ \\
\gafp^6+\gbfp^6 & = & \glfp^6 & - & 6\glfp^4 & + & \ \ 9\glfp^2 & -
&\ \
2  &  &  \\ \ \\
\gafp^7+\gbfp^7 & = & \glfp^7 & - & 7\glfp^5 & + & 14\glfp^3 & - & \
\
7\glfp &  & \\ \ \\
\gafp^8+\gbfp^8 & = & \glfp^8 & - & 8\glfp^6 & + & 20\glfp^4 & - & 16\glfp^2 & + & 2. \\
                    \end{array} \eea Writing $\gafp^m +
\gbfp^m$ as a polynomial in $\glfp$, we find that \be \gafp^m +
\gbfp^m \ = \ \sum_{r=0\atop r\equiv m \bmod 2}^m c_{m,r} \glfp^r,
\ee where the $c_{m,r}$ are our coefficients from Definition
\ref{def:formscmr}. A key ingredient in the proof is noting that
 \ben  \item $c_{2k,2\ell} = c_{2k-1,2\ell-1} -
c_{2k-2,2\ell}$ if $\ell \in \{1,\dots,k-1\}$ and $k \ge 2$; \item
$c_{2k+1,2\ell+1} = c_{2k,2\ell} - c_{2k-1,2\ell+1}$ if $\ell < k$.
\een

We briefly describe the application of our identity, ignoring the book-keeping needed to deal with $m \le 2$. From the explicit formula \eqref{eq:explicitformula}, we see we must understand sums such as \be \sum_p
\sum_{m=3}^\infty \frac1{W_R(\F)} \sum_{f\in\F}w_R(f)\frac{\alpha_f(p)^m + \beta_f(p)^m}{p^{m/2}} \frac{\log
p}{\log R}\ \hphipr, \ee where $\F$ is a family of cuspidal newforms and $W_R(\F) = \sum_{f\in \F} w_R(f)$ (a simple Taylor series shows there is negligible contribution in replacing $\hphi(m\log p /\log R)$ with $\hphi(\log p/\log R)$). As the sums of powers of the Satake parameters are polynomials in $\glfp$, we may rewrite this as \be  \sum_p
\sum_{m=3}^\infty \sum_{r=0 \atop r \equiv m \bmod 2}^m
\frac{c_{m,r} A_{r,\F}(p)}{p^{m/2}} \frac{\log p}{\log R}\
\hphipr,  \ee where $A_{r,\F}(p)$ is the $r$\textsuperscript{th} moment of $\glfp$ in the family $\F$: \bea A_{r,\F}(p) & \ = \ & \fwf\sum_{f \in
\F \atop f \in S(p)} w_R(f)\glfp^r. \eea We interchange the $m$ and $r$ sums (which is straightforward for $p \ge 11$, and follows by Abel summation for $p \le 7$) and then apply our polylogarithm identity (Theorem \ref{thm:lincombpolylog}) to rewrite the sum as \be \sum_{p} \sum_{r = 0}^\infty
\frac{A_{r,\F}(p)p^{r/2}(p-1)\log p}{(p+1)^{r+1}\log R} \ \hphipr. \ee

For many families we either know or conjecture a distribution for
the (weighted) Fourier coefficients. If this were the case, then we
could replace the $A_{r,\F}(p)$ with the $r$\textsuperscript{th}
moment. In many applications (for example, using the Petersson
formula for families of cuspidal newforms of fixed weight and
square-free level tending to infinity) we know the moments up to a
negligible correction (the distribution is often known or conjectured to be Sato-Tate, unless we are looking at families of elliptic curves with complex multiplication, where the distribution is known and slightly more complicated). Simple algebra yields

\begin{lem} Assume for $r \ge 3$ that \be \twocase{A_{r,\F}(p) \ = \ }{M_\ell +
O\left(\frac1{\log^2 R}\right)}{if $r = 2\ell$}{0}{otherwise,}\ee
and that there is a nice function $g_M$ such that \be g_M(x) \ = \
M_2 x^2 + M_3 x^3 + \cdots \ = \  \sum_{\ell=2}^\infty M_\ell\
x^\ell. \ee Then the contribution from the $r \ge 3$ terms in the explicit formula is \be -\frac{2\hphi(0)}{\log R} \sum_p
g_M\left( \frac{p}{(p+1)^2}\right) \cdot \frac{(p-1)\log p}{p+1} +
O\left(\frac1{\log^3 R}\right). \ee
\end{lem}

Thus we can use our polylogarithm identity to rewrite the sums arising in the explicit formula in a very compact way which emphasizes properties of the known or conjectured distribution of the Fourier coefficients. One application of this is in analyzing the behavior of the zeros of $L$-functions near the central point. Many investigations have shown that, for numerous families, as the conductors tend to infinity the behavior of these zeros is the same as the $N\to\infty$ scaling limit of eigenvalues near $1$ of subgroups of $U(N)$.

Most of these studies only examine the main term, showing agreement in the limit with random matrix theory (the scaling limits of eigenvalues of $U(N)$). In particular, all one-parameter families of elliptic curves over $\Q(T)$ with the same rank and same limiting distribution of signs of functional equation have the same main term for the behavior of their zeros. What is unsatisfying about this is that the arithmetic of the families is not seen; this is remedied, however, by studying the lower order terms in the $1$-level density. There we \emph{do} break the universality and see arithmetic dependent terms. In particular, our formula shows that we have different answers for families of elliptic curves with and without complex multiplication (as these two cases have different densities for the Fourier coefficients).

These lower order differences, which reflect the arithmetic structure of the family, are quite important. While the behavior of many properties of zeros of $L$-functions of height $T$ are well-modeled by the $N\to\infty$ scaling limits of eigenvalues of a classical compact group, better agreement (taking into account lower order terms) is given by studying matrices of size $N = (\log T)/2\pi$ (see \cite{KeSn1,KeSn2,KeSn3}). Recently it has been observed that even better agreement is obtained by replacing $N$ with $N_{\rm eff}$, where $N_{\rm eff}$ is chosen so that the main and first lower order terms match (see \cite{BBLM,DHKMS}). Thus one consequence of our work is in deriving a tractable formula to identify the lower order correction terms, which results in an improved model for the behavior of the zeros.


\bigskip

\end{document}